\documentclass{article}

\author{Dixy Msapato}
\usepackage{amsmath,enumerate,amsthm,amsfonts,amssymb,latexsym,bbm}
\usepackage[super]{nth}
\usepackage[a4paper]{geometry}
\usepackage[hidelinks]{hyperref}
\usepackage[UKenglish]{babel}
\usepackage{float}
\usepackage{dirtytalk}
\usepackage[mathscr]{euscript} %fancy mathcal
\usepackage{tikz}
\usepackage{tikz-cd}
\usepackage{tikz-qtree}
\usepackage{mathtools}
\usepackage{xcolor}

\newcommand{\myfrac}%  % macro with LaTeX-style syntax
    [2]{\begin{array}{@{}c@{}}#1 \\[-0.75ex]#2\end{array}}

\providecommand{\Ext}{\mathop{\rm Ext}\nolimits}

\newdimen\R
\newdimen\T
\newcommand\blfootnote[1]{%
  \begingroup
  \renewcommand\thefootnote{}\footnote{#1}%
  \addtocounter{footnote}{-1}%
  \endgroup
}

\theoremstyle{definition}
\newtheorem{theorem}{Theorem}[section]
\newtheorem{definition}[theorem]{Definition}
\newtheorem{proposition}[theorem]{Proposition}
\newtheorem{remark}[theorem]{Remark}

\newtheorem{lemma}[theorem]{Lemma}

\def\d{\delta}

\title{A characterisation of extriangulated categories with triangulated structure}
\date{}

\begin{document}

\maketitle
\begin{abstract}
We give a characterisation of the extriangulated categories which admit the structure of a triangulated category. We show that these are the extriangulated categories where for every object $X$ in the extriangulated category, the morphism $0 \rightarrow X$ is a deflation and the morphism $X \rightarrow 0$ is an inflation.
\end{abstract}


%\tableofcontents
\thanks{}
\blfootnote{Address: School of Mathematics, University of Leeds, Leeds, LS2 9JT, United Kingdom. \\
	Contact: mmdmm@leeds.ac.uk \\
	\thanks{This research was supported by an EPSRC Doctoral Training Partnership (reference EP/R513258/1) through the University of Leeds. The author also wishes to thank their supervisor, Bethany Marsh for their continuous support and invaluable insight. The author also expresses gratitude to Amit Shah for their helpful discussions. }}
\section{Introduction.}
Extriangulated categories were introduced by Nakaoka and Palu in \cite{NakaokaPalu} as a simultaneous generalisation of exact categories and triangulated categories in the context of the study of cortosion pairs. The known classes of examples of extriangulated categories include exact categories and extension-closed subcategories of triangulated categories; see \cite[Example 2.13, Remark 2.18, Proposition 3.22(1)]{NakaokaPalu}. There are also examples of extriangulated categories which are neither exact nor triangulated; see for example, \cite[Proposition 3.30]{NakaokaPalu}, \cite[Example 4.14 and Corollary 4.12]{ZhouZhu}.

The extriangulated categories which also have an exact or triangulated structure have been characterised. The extriangulated categories which are exact categories are those where every inflation is a monomorphism and every deflation is an epimorphism; see \cite[Corollary 3.18]{NakaokaPalu}. The extriangulated categories which are triangulated are those where the $\Ext^{1}$ bifunctor $\mathbb{E}$ is such that $\mathbb{E}(-,-) = \text{Hom}(-,\Sigma -)$ for some auto-equivalence $\Sigma$ on the category; see \cite[Proposition 3.22(2)]{NakaokaPalu}. In this paper, we offer another characterisation of the extriangulated categories which are triangulated, which is as follows. An extriangulated category has the structure of a triangulated category if and only for all objects $X$ in the category, the morphism $0 \rightarrow X$ is a deflation and the morphism $X \rightarrow 0$ is an inflation. See Theorem \ref{mainTheorem} of this paper.

This paper is organised as follows: in \S 2, we recall the necessary theory of extriangulated categories. In \S 3, we state and prove our main theorem, a characterisation of the extriangulated categories which have a structure of a triangulated category. 

\section{Extriangulated categories.}
In this section, we will recall from \cite{NakaokaPalu} the basic theory of extriangulated categories needed for this paper. Through out this subsection, $\mathscr{C}$ will be an additive category equipped with a biadditive functor $\mathbb{E} \colon \mathscr{C}^{\text{op}} \times \mathscr{C} \rightarrow Ab$, where $Ab$ is the category of abelian groups. 

\begin{definition}\cite[Definition 2.1]{NakaokaPalu}.
Let $A,C$ be objects of $\mathscr{C}$. An element $\d \in \mathbb{E}(C,A)$ is called an $\mathbb{E}$-\textit{extension}. Formally an $\mathbb{E}$-extension $\d \in \mathbb{E}(C,A)$ is a triple $(A,\d,C)$. 

Since $\mathbb{E}$ is a bifunctor, for any morphisms $a \in \mathscr{C}(A,A^{\prime})$ and $c \in \mathscr{C}(C^{\prime},C)$, we have the following $\mathbb{E}$-extensions:
$$a_{*}\d := \mathbb{E}(C,a)(\d) \in \mathbb{E}(C,A^{\prime}),$$
$$c^{*}\d := \mathbb{E}(c^{\text{op}},A)(\d) \in \mathbb{E}(C^{\prime},A) \text{ and }$$
$$c^{*}a_{*}\d = a_{*}c^{*}\d := \mathbb{E}(c^{\text{op}},a)(\d) \in \mathbb{E}(C^{\prime},A^{\prime}).$$
\end{definition}
For the rest of the paper, we will abuse notation by writing $\mathbb{E}(c,-)$ instead of $\mathbb{E}(c^{\text{op}},-).$

\begin{definition}\cite[Definition 2.3]{NakaokaPalu}.
Let $(A,\d,C)$ and $(A^{\prime}, \d^{\prime}, C^{\prime})$ be any pair of $\mathbb{E}$-extensions. A morphism of $\mathbb{E}$-extensions $(a,c) \colon \d \rightarrow \d^{\prime}$ is a pair of morphisms $a \in \mathscr{C}(A,A^{\prime})$ and $c \in \mathscr{C}(C,C^{\prime})$ such that: $$ a_{*}\d = c^{*}\d^{\prime}.$$ 
\end{definition}

\begin{lemma}\label{*operation}\cite[Remark 2.4]{NakaokaPalu}.
Let $(A,\d,C)$ be an $\mathbb{E}$-extension. Then we have the following.
\begin{enumerate}
\item Any morphism $a \in \mathscr{C}(A,A^{\prime})$ induces the following morphism of $\mathbb{E}$-extensions,
$$(a,1_C) \colon \d \rightarrow a_{*}\d.$$
\item Any morphism $c \in \mathscr{C}(C^{\prime},C)$ induces the following morphism of $\mathbb{E}$-extensions,
$$(1_A,c) \colon c^{*}\d \rightarrow \d.$$
\end{enumerate}
\end{lemma}

\begin{definition}\cite[Definition 2.5]{NakaokaPalu}.
For any objects $A,C$ in $\mathscr{C}$, the zero element $0 \in \mathbb{E}(C,A)$ is called the \textit{split $\mathbb{E}$-extension}. 
\end{definition}

\begin{definition}\cite[Definition 2.6]{NakaokaPalu}.
Let $\d \in \mathbb{E}(C,A)$ and $\d^{\prime} \in \mathbb{E}(C^{\prime}, A^{\prime})$ be any pair of $\mathbb{E}$-extensions. Let $i_{C} \colon C \rightarrow C \oplus C^{\prime}$ and $i_{C^{\prime}} \colon C^{\prime} \rightarrow C \oplus C^{\prime}$ be the canonical inclusion maps. Let $p_{A} \colon A \oplus A^{\prime} \rightarrow A$,  and $p_{A^{\prime}} \colon A \oplus A^{\prime} \rightarrow A^{\prime}$ be the canonical projection maps. By the biadditivity of $\mathbb{E}$ we have the following isomorphism.
$$ \mathbb{E}(C \oplus C^{\prime}, A \oplus A^{\prime}) \cong \mathbb{E}(C,A) \oplus \mathbb{E}(C,A^{\prime}) \oplus \mathbb{E}(C^{\prime}, A) \oplus \mathbb{E}(C^{\prime}, A^{\prime})$$ 

Let $\d \oplus \d^{\prime} \in \mathbb{E}(C \oplus C^{\prime}, A \oplus A^{\prime})$ be the element corresponding to $(\d,0,0,\d^{\prime})$ via the above isomorphism. 
If $A = A^{\prime}$ and $C = C^{\prime}$, then the sum $\d + \d^{\prime} \in \mathbb{E}(C,A)$ is obtained by 
$$ \d + \d^{\prime} = \mathbb{E}(\Delta_{C},\nabla_{A})(\d \oplus \d^{\prime}),$$
where $\Delta_{C} = \begin{pmatrix}
1 \\ 1
\end{pmatrix} : C \rightarrow C \oplus C$, and $\nabla_{A} = \begin{pmatrix}
1,1
\end{pmatrix}: A \oplus A \rightarrow A.$
\end{definition}

\begin{definition}\label{equivReln}\cite[Definition 2.7]{NakaokaPalu}.
Let $A,C$ be a pair of objects in $\mathscr{C}$. Two sequences of morphisms $A \overset{x}{\longrightarrow} B \overset{y}{\longrightarrow} C$ and $A \overset{x^{\prime}}{\longrightarrow} B^{\prime} \overset{y^{\prime}}{\longrightarrow} C$ in $\mathscr{C}$ are said to be equivalent if there exists an isomorphism $b \in \mathscr{C}(B,B^{\prime})$ such that the following diagram commutes. 
\begin{center}
\begin{tikzcd}
A \arrow[r, "x"] \arrow[d, equal]
& B \arrow[d, "b"] \arrow[r,"y"]
& C \arrow[d, equal] \\
 A \arrow[r,"x^{\prime}"]
&B^{\prime} \arrow[r,"y^{\prime}"]
&C \end{tikzcd}
\end{center}
We denote the equivalence class of a sequence  $A \overset{x}{\longrightarrow} B \overset{y}{\longrightarrow} C$ by 
$[A \overset{x}{\longrightarrow} B \overset{y}{\longrightarrow} C]$.
\end{definition}

\begin{definition}\cite[Definition 2.8]{NakaokaPalu}.
Let $A,B,C,A^{\prime}, B^{\prime}, C^{\prime}$ be objects in the category $\mathscr{C}$. 
\begin{enumerate}
\item We denote by $0$ the equivalence class $[A \overset{\big[\begin{smallmatrix}
1_A\\
0
\end{smallmatrix}\big]}{\longrightarrow} A \oplus C \overset{[\begin{smallmatrix} 0 & 1_C \end{smallmatrix}]}{\longrightarrow} C]$.

\item For any two equivalence classes  $[A \overset{x}{\longrightarrow} B \overset{y}{\longrightarrow} C]$ and $[A^{\prime} \overset{x^{\prime}}{\longrightarrow} B^{\prime} \overset{y^{\prime}}{\longrightarrow} C^{\prime}]$, we denote by $[A \overset{x}{\longrightarrow} B \overset{y}{\longrightarrow} C] \oplus [A^{\prime} \overset{x^{\prime}}{\longrightarrow} B^{\prime} \overset{y^{\prime}}{\longrightarrow} C^{\prime}]$ the equivalence class
$[A \oplus A^{\prime} \overset{x\oplus x^{\prime}}{\longrightarrow} B \oplus B^{\prime} \overset{y \oplus y^{\prime}}{\longrightarrow} C \oplus C^{\prime}].$
\end{enumerate}
\end{definition}

\begin{definition}\cite[Definition 2.9]{NakaokaPalu}.
Let $\mathfrak{s}$ be a correspondence associating an equivalence class $\mathfrak{s}(\d) = [A \overset{x}{\longrightarrow} B \overset{y}{\longrightarrow} C]$ to any $\mathbb{E}$-extension $\d \in \mathbb{E}(C,A)$.  We say that $\mathfrak{s}$ is a \textit{realisation} of $\mathbb{E}$ if the following condition $(\circ)$ holds. 

$(\circ)$ Let $\d \in \mathbb{E}(C,A)$ and $\d^{\prime} \in \mathbb{E}(C^{\prime},A^{\prime})$ be $\mathbb{E}$-extensions with $\mathfrak{s}(\d) = [A \overset{x}{\longrightarrow} B \overset{y}{\longrightarrow} C]$ and $\mathfrak{s}(\d^{\prime}) = [A^{\prime} \overset{x^{\prime}}{\longrightarrow} B^{\prime} \overset{y^{\prime}}{\longrightarrow} C^{\prime}]$. Then for any morphism of $\mathbb{E}$-extensions $(a,c) \colon \d \rightarrow \d^{\prime}$, there exists a morphism $b \in \mathscr{C}(B,B^{\prime})$ such that the following diagram commutes. 
\begin{center}
\begin{tikzcd}
A \arrow[r, "x"] \arrow[d, "a"]
& B \arrow[d, dashed, "b"] \arrow[r,"y"]
& C \arrow[d, "c"] \\
 A^{\prime} \arrow[r,"x^{\prime}"]
&B^{\prime} \arrow[r,"y^{\prime}"]
&C^{\prime} \end{tikzcd}
\end{center}
In this situation, we say that the triple of morphisms $(a,b,c)$ \textit{realises} $(a,c)$. Moreover, for any $\d \in \mathbb{E}(C,A)$, we say that the sequence $A \overset{x}{\longrightarrow} B \overset{y}{\longrightarrow} C$ realises $\d$ if $\mathfrak{s}(\d)=[A \overset{x}{\longrightarrow} B \overset{y}{\longrightarrow} C]$.
\end{definition}

\begin{definition}\label{additiveRealisation}\cite[Definition 2.16]{HassounShah}.
A realisation $\mathfrak{s}$ is said to be an \textit{additive realisation} if the following conditions are satisfied,
\begin{enumerate}
\item For any objects $A,C$ in $\mathscr{C}$, a split $\mathbb{E}$-extension $ 0 \in \mathbb{E}(C,A)$ satisfies $$\mathfrak{s}(0)=0.$$
\item Let $\d \in \mathbb{E}(C,A)$ and $\d^{\prime} \in \mathbb{E}(C^{\prime},A^{\prime})$ be $\mathbb{E}$-extensions with $\mathfrak{s}(\d) = [A \overset{x}{\longrightarrow} B \overset{y}{\longrightarrow} C]$ and $\mathfrak{s}(\d^{\prime}) = [A^{\prime} \overset{x^{\prime}}{\longrightarrow} B^{\prime} \overset{y^{\prime}}{\longrightarrow} C^{\prime}]$. Let $i_{A} \colon A \rightarrow A \oplus A^{\prime}$ and $i_{A^{\prime}} \colon A^{\prime} \rightarrow A \oplus A^{\prime}$ be the canonical inclusions. Let $p_{C} \colon C \oplus C^{\prime} \rightarrow C$ and $p_{C^{\prime}} \colon C \oplus C^{\prime} \rightarrow C^{\prime}$ be the canonical projections, then the $\mathbb{E}$-extension
$$(i_{A})_{*}(p_{C})^{*}\d + (i_{A^{\prime}})_{*}(p_{C^{\prime}})^{*}\d^{\prime} \in \mathbb{E}(C \oplus C^{\prime}, A \oplus A^{\prime})$$ is realised by the direct sum
\begin{center}
\begin{tikzcd}
A \oplus A^{\prime} \arrow[r,"x \oplus x^{\prime}"]
& B \oplus B^{\prime} \arrow[r,"y \oplus y^{\prime}"]
& C \oplus C^{\prime}.
\end{tikzcd}
\end{center}
\end{enumerate} 
\end{definition}

\begin{remark}
The reader familiar with the theory of extriangulated categories or the seminal paper by Nakaoka and Palu \cite{NakaokaPalu} will notice that the definition of an additive realisation above is not stated in the usual way as in \cite[Definition 2.10]{NakaokaPalu}. In \cite[Definition 2.10]{NakaokaPalu}, the second part of the above definition is given as the statement $\mathfrak{s}(\d \oplus \d^{\prime}) = \mathfrak{s}(\d) \oplus \mathfrak{s}(\d^{\prime})$. For our main purpose in this paper, the restatement of the definition of an additive realisation as in \cite[Definition 2.16]{HassounShah} turns out to be more convenient. 
\end{remark}
We are now in a position to define an extriangulated category.

\begin{definition}\label{DefExtriang}\cite[Definition 2.12]{NakaokaPalu}. Let $\mathscr{C}$ be an additive category. An \textit{extriangulated category} is a triple $(\mathscr{C},\mathbb{E},\mathfrak{s})$ satisfying the following axioms.

\begin{itemize}

\item[] (ET1) The functor $\mathbb{E} \colon \mathscr{C}^{\text{op}} \times \mathscr{C} \rightarrow Ab \text{ is a biadditive functor}$. 

\item[] (ET2) The correspondence $\mathfrak{s}$ is an additive realisation of $\mathbb{E}$. 

\item[] (ET3) Let $\d \in \mathbb{E}(C,A)$ and $\d^{\prime} \in \mathbb{E}(C^{\prime},A^{\prime})$ be any pair of $\mathbb{E}$-extensions realised by the sequences $A \overset{x}{\longrightarrow} B \overset{y}{\longrightarrow} C$ and $A^{\prime} \overset{x^{\prime}}{\longrightarrow} B^{\prime} \overset{y^{\prime}}{\longrightarrow} C^{\prime}$ respectively. Then for any commutative diagram 
\begin{center}
\begin{tikzcd}
A \arrow[r, "x"] \arrow[d, "a"]
& B \arrow[d, "b"] \arrow[r,"y"]
& C \\
 A^{\prime} \arrow[r,"x^{\prime}"]
&B^{\prime} \arrow[r,"y^{\prime}"]
&C^{\prime} \end{tikzcd}
\end{center}
there exists a morphism $c \in \mathscr{C}(C,C^{\prime})$ such that $(a,c) \colon \d \rightarrow \d^{\prime}$ is a morphism of $\mathbb{E}$-extensions and the triple $(a,b,c)$ realises $(a,c)$.

\item[] $\text{(ET3})^{\text{op}}$ The dual of (ET3).

\item[] (ET4) Let $\d \in \mathbb{E}(D,A)$ and $\d^{\prime} \in \mathbb{E}(F,B)$ be any pair of $\mathbb{E}$-extensions realised by the sequences, $A \overset{f}{\longrightarrow} B \overset{f^{\prime}}{\longrightarrow} D$ and $B\overset{g}{\longrightarrow} C \overset{g^{\prime}}{\longrightarrow} F$ respectively. Then there exists an object $E$ in $\mathscr{C}$, a commutative diagram
\begin{center}
\begin{tikzcd}
A \arrow[r, "f"] \arrow[d, equal]
& B \arrow[d, "g"] \arrow[r,"f^{\prime}"]
& D \arrow[d,"d"] \\
 A \arrow[r,"h"]
&C \arrow[r,"h^{\prime}"] \arrow[d,"g^{\prime}"]
&E \arrow[d,"e"] \\
 & F \arrow[r,equal] & F
\end{tikzcd}
\end{center}
in $\mathscr{C}$ and an $\mathbb{E}$-extension $\d^{\prime \prime} \in \mathbb{E}(E,A)$ realised by the sequence $A \overset{h}{\longrightarrow} C \overset{h^{\prime}}{\longrightarrow} E$, such that the following compatibilities are satisfied;
\begin{enumerate}[(i)]
\item $\mathfrak{s}((f^{\prime})_{*}\d^{\prime}) = [D \overset{d}{\longrightarrow} E \overset{e}{\longrightarrow} F].$

\item $d^{*} \d^{\prime \prime} = \d.$

\item $f_{*}\d^{\prime \prime} = e^{*}\d^{\prime}.$
\end{enumerate}

\item[] $\text{(ET4)}^{\text{op}}$ The dual of (ET4). 
\end{itemize}
In this case, we call $\mathfrak{s}$ an $\mathbb{E}$-\textit{triangulation} of $\mathscr{C}$.
\end{definition}
There are many examples of extriangulated categories. They include exact categories, triangulated categories and extension-closed subcategories of triangulated subcategories. There are also extriangulated categories which are neither exact nor triangulated; for examples see, \cite[Proposition 3.30]{NakaokaPalu}, \cite[Example 4.14 and Corollary 4.12]{ZhouZhu}.

We will conclude this section by introducing some useful terminology from \cite{NakaokaPalu} and stating results about extriangulated categories which will be helpful for the rest of the paper. 

\begin{definition}\cite[Definition 2.5, Definition 3.9]{NakaokaPalu}.
Let $(\mathscr{C},\mathbb{E},\mathfrak{s})$ be a triple satisfying (ET1) and (ET2). 
\begin{enumerate}
\item A sequence $A \overset{x}{\longrightarrow} B \overset{y}{\longrightarrow} C$ is called an \textit{$\mathbb{E}$-conflation} if it realises some $\mathbb{E}$-extension $\d \in \mathbb{E}(C,A)$. When there is no risk of confusion, we will refer to $\mathbb{E}$-conflations simply as \textit{conflations}.
\item A morphism $f \in \mathscr{C}(A,B)$ is called an \textit{$\mathbb{E}$-inflation} if it admits some conflation $A \overset{f}{\longrightarrow} B \longrightarrow C$. In this case, we call $C$ \textit{a cone of f} and denote it by Cone$(f)$. When there is no risk of confusion, we will refer to $\mathbb{E}$-inflations simply as \textit{inflations.}
\item A morphism $g \in \mathscr{C}(B,C)$ is called an \textit{$\mathbb{E}$-deflation} if it admits some conflation $A \longrightarrow B \overset{g}{\longrightarrow} C$. In this case, we call $A$ \textit{a cocone of g} and denote it by Cocone$(g)$. When there is no risk of confusion, we will refer to $\mathbb{E}$-deflations simply as \textit{deflations.}
\end{enumerate}
The objects Cone$(f)$ and Cocone$(g)$ are unique up to isomorphism; see \cite[Remark 3.10]{NakaokaPalu}.
\end{definition}
The terminology of conflations, inflations and deflations is also used in the context of exact categories and triangulated categories analogously. 

\begin{definition}\cite[Definition 2.19]{NakaokaPalu}. Let $(\mathscr{C},\mathbb{E},\mathfrak{s})$ be a triple satisfying (ET1) and (ET2). 
\begin{enumerate}
\item If a conflation $A \overset{x}{\longrightarrow} B \overset{y}{\longrightarrow} C$ realises $\d \in \mathbb{E}(C,A)$, we call the pair $(A \overset{x}{\longrightarrow} B \overset{y}{\longrightarrow} C, \d)$ an $\mathbb{E}$-\textit{triangle} and denote it by the following diagram.
\begin{center}
\begin{tikzcd}
A \arrow[r, "x"]
& B \arrow[r,"y"]
& C \arrow[r, dashed, "\d"] & \text{}
\end{tikzcd} 
\end{center}

\item Let \begin{tikzcd} 
A \arrow[r, "x"]
& B \arrow[r,"y"]
& C \arrow[r, dashed, "\d"] &\text{}
\end{tikzcd}and
\begin{tikzcd} 
A^{\prime} \arrow[r, "x^{\prime}"]
& B^{\prime} \arrow[r,"y^{\prime}"]
& C^{\prime} \arrow[r, dashed, "\d^{\prime}"] &\text{}
\end{tikzcd}
be any pair of $\mathbb{E}$-triangles. If a triple $(a,b,c)$ realises $(a,c) \colon \d \rightarrow \d^{\prime}$ we write it as in the following commutative diagram and call $(a,b,c)$ a morphism of $\mathbb{E}$-triangles. 
\begin{center}
\begin{tikzcd}
A \arrow[r, "x"] \arrow[d,"a"]
& B \arrow[r,"y"] \arrow[d,"b"]
& C \arrow[r, dashed, "\d"] \arrow[d,"c"] & \text{} \\
A^{\prime} \arrow[r, "x^{\prime}"]
& B^{\prime} \arrow[r,"y^{\prime}"]
& C^{\prime} \arrow[r, dashed, "\d^{\prime}"] &\text{}
\end{tikzcd}
\end{center}

\end{enumerate}
\end{definition}

\begin{lemma}\label{2out3iso}\cite[Corollary 3.6]{NakaokaPalu}. 
Let $(\mathscr{C},\mathbb{E},\mathfrak{s})$ be a triple satisfying (ET1), (ET2), (ET3) and $\text{(ET3)}^{\text{op}}$. Let $(a,b,c)$ be a morphism of $\mathbb{E}$-triangles. If any two of $a,b,c$ are isomorphisms, then so is the third. 
\end{lemma}

\begin{proposition}\label{longExact}\cite[Corollary 3.12]{NakaokaPalu}.
Let $(\mathscr{C},\mathbb{E},\mathfrak{s})$ be an extriangulated category. For any $\mathbb{E}$-triangle \begin{tikzcd} 
A \arrow[r, "x"]
& B \arrow[r,"y"]
& C \arrow[r, dashed, "\d"] &\text{,}
\end{tikzcd}the following sequences of natural transformations are exact. 

\begin{center}
\begin{tikzcd}
\mathscr{C}(C,-) \arrow[r, Rightarrow, "\mathscr{C}(y{,}-)"]
& \mathscr{C}(B,-) \arrow[r, Rightarrow, "\mathscr{C}(x{,}-)"]
& \mathscr{C}(A,-) \arrow[r, Rightarrow, "\d^{\#}"]
& \mathbb{E}(C,-) \arrow[r, Rightarrow, "\mathbb{E}(y{,}-)"]
& \mathbb{E}(B,-) \arrow[r, Rightarrow, "\mathbb{E}(x{,}-)"]
& \mathbb{E}(A,-)
\end{tikzcd}
\end{center}

\begin{center}
\begin{tikzcd}
\mathscr{C}(-,A) \arrow[r, Rightarrow, "\mathscr{C}(-{,}x)"]
& \mathscr{C}(-,B) \arrow[r, Rightarrow, "\mathscr{C}(-{,}y)"]
& \mathscr{C}(-,C) \arrow[r, Rightarrow, "\d_{\#}"]
& \mathbb{E}(-,A) \arrow[r, Rightarrow, "\mathbb{E}(-{,}x)"]
& \mathbb{E}(-,B) \arrow[r, Rightarrow, "\mathbb{E}(-{,}y)"]
& \mathbb{E}(-,C)
\end{tikzcd}
\end{center}
The natural transformations $\d^{\#}$ and $\d_{\#}$ are defined as follows. Given any object $X$ in $\mathscr{C}$, we have that
\begin{enumerate}
\item $(\d^{\#})_{X} \colon \mathscr{C}(A,X) \rightarrow \mathbb{E}(C,X) \text{ ; } g \mapsto g_{*}\d,$

\item $(\d_{\#})_{X} \colon \mathscr{C}(X,C) \rightarrow \mathbb{E}(X,A) \text{ ; } f \mapsto f^{*}\d.$
\end{enumerate}

The exactness of the first sequence of natural transformations is taken to mean that for any object $X$ in $\mathscr{C}$, the sequence
\begin{center}
\begin{tikzcd}
\mathscr{C}(C,X) \arrow[r, "\mathscr{C}(y{,}X)"]
& \mathscr{C}(B,X) \arrow[r, "\mathscr{C}(x{,}X)"]
& \mathscr{C}(A,X) \arrow[r, "\d^{\#}_{X}"]
& \mathbb{E}(C,X) \arrow[r, "\mathbb{E}(y{,}X)"]
& \mathbb{E}(B,X) \arrow[r,"\mathbb{E}(x{,}X)"]
& \mathbb{E}(A,X)
\end{tikzcd}
\end{center}
is exact in $Ab$ and likewise for the second sequence.
\end{proposition}

\section{A characterisation of triangulated extriangulated categories.}
This section is dedicated to giving a new characterisation of the extriangulated categories which are triangulated. In \cite{NakaokaPalu}, the authors characterise the extriangulated categories which are exact and the ones which are triangulated. They characterise the extriangulated categories which are triangulated as those where $\mathbb{E} = \mathscr{C}(-,\Sigma -)$ for some auto-equivalence $\Sigma \colon \mathscr{C} \rightarrow \mathscr{C}$. They also characterise the extriangulated categories which have the structure of an exact category as those where every inflation is a monomorphism and every deflation is an epimorphism. The forward direction is by the definition of exact categories, since inflations are kernels and deflations are cokernels; see \cite[Definition 2.1]{buhler}. The other direction was shown by Nakaoka and Palu, see \cite[Corollary 3.18]{NakaokaPalu}.  In this section, we offer a characterisation of the extriangulated categories which are triangulated that is in the spirit of the characterisation of the extriangulated categories which are exact. 

For the benefit of the reader, we start by recalling the definition of a triangulated category and some terminology. Our reference for triangulated categories is \cite[\S 1]{Happel}. 

Let $\mathscr{A}$ be an additive category and $T$ an auto-equivalence of $\mathscr{A}$. A sextuple $(X,Y,Z,u,v,w)$ in $\mathscr{A}$ is a sequence of morphisms of the following form
\begin{center}
\begin{tikzcd}
X \arrow[r, "u"]
& Y \arrow[r,"v"]
& Z \arrow[r,"w"] 
& T(X)
\end{tikzcd} 
\end{center}
A morphism of sextuples from $(X,Y,Z,u,v,w)$ to $(X^{\prime}, Y^{\prime}, Z^{\prime}, u^{\prime}, v^{\prime}, w^{\prime})$ is a triple $(f,g,h)$ of morphisms such that the following diagram commutes. 
\begin{center}
\begin{tikzcd}
X \arrow[r,"u"] \arrow[d,"f"]
&
Y \arrow[r, "v"] \arrow[d,"g"]
&
Z \arrow[r, "w"] \arrow[d, "h"]
& 
TX \arrow[d, "Tf"]
\\
X^{\prime} \arrow[r,"u^{\prime}"]
&
Y^{\prime} \arrow[r, "v^{\prime}"]
&
Z^{\prime} \arrow[r, "w^{\prime}"]
& 
TX^{\prime}
\end{tikzcd}
\end{center}
If $f,g $ and $h$ are isomorphisms in $\mathscr{A}$, we say that the triple $(f,g,h)$ is an isomorphism. 

\begin{definition}\cite[\S 1]{Happel} Let $\mathscr{A}$ be an additive category with auto-equivalence $T$. Let $\mathcal{T}$ be a set of sextuples of $\mathscr{A}$. The triple $(\mathscr{A},T, \mathcal{T})$ is called a triangulated category if the following axioms hold. In this case, the elements of $\mathcal{T}$ are then called \textit{distinguished triangles}.
\begin{itemize}
\item[] (TR1) Every sextuple isomorphic to a distinguished triangle is again a distinguished triangle. Every morphism $u \colon X \rightarrow Y$ in $\mathscr{A}$ can be embedded into a distinguished triangle $(X,Y,Z,u,v,w)$. The sextuple $(X,X,0,1_X,0,0)$ is a distinguished triangle. 
\item[] (TR2) If $(X,Y,Z,u,v,w)$ is a distinguished triangle, then $(Y,Z,TX,v,w,-Tu)$ is also a distinguished triangle. 
\item[] (TR3) Given two distinguished triangles $(X,Y,Z,u,v,w), (X^{\prime}, Y^{\prime}, Z^{\prime}, u^{\prime}, v^{\prime}, w^{\prime})$ and two morphisms $f \colon X \rightarrow X^{\prime}$, $g \colon Y \rightarrow Y^{\prime}$ such that $u^{\prime}f = gu$. Then there exists a morphism $h \colon Z \rightarrow Z^{\prime}$ such that $(f,g,h)$ is a morphism of distinguished from the first distinguished triangle to the second.
\item[] (TR4) Consider the distinguished triangles $$(X,Y,Z^{\prime},u,i,i^{\prime}), (Y,Z,X^{\prime},v,j,j^{\prime}) \text{ and } (X,Z,Y^{\prime}, u \circ v, k,k^{\prime}).$$ Then there exists morphisms $f \colon Z^{\prime} \rightarrow Y^{\prime}$, $g \colon Y^{\prime} \rightarrow X^{\prime}$ such that the following diagram commutes and the second column is a distinguished triangle. Moreover, we have that $u[1]k^{\prime} = j^{\prime}g$.
\begin{center}
\begin{tikzcd}
X \arrow[r, "u"] \arrow[d,equal]
& Y \arrow[r, "i"] \arrow[d,"v"]
&Z^{\prime} \arrow[r,"i^{\prime}"] \arrow[d,"f"]
&TX \arrow[d, equal]
\\
X \arrow[r, "v\circ u"]
& Z \arrow[r, "k"] \arrow[d,"j"]
& Y^{\prime} \arrow[r,"k^{\prime}"] \arrow[d,"g"]
&TX
\\
\text{}
& X^{\prime} \arrow[r,equal] \arrow[d, "j^{\prime}"]
& X^{\prime} \arrow[d,"Ti \circ j^{\prime}"]
\\
\text{}
& TY \arrow[r,"Ti"]
&
TZ

\end{tikzcd}
\end{center}

\end{itemize}

\end{definition}

\begin{definition} Let $(\mathscr{C},\mathbb{E},\mathfrak{s})$ be an extriangulated category. Suppose $\mathscr{C}$ has the structure of a triangulated category $(\mathscr{C},T,\mathcal{T})$. We say that this triangulated structure is $\mathbb{E}$-\textit{compatible} if and only if for each distinguished triangle 
\begin{center}
\begin{tikzcd}
X \arrow[r, "u"]
& Y \arrow[r,"v"]
& Z \arrow[r,"w"] 
& T(X)
\end{tikzcd} 
\end{center}

we have that
\begin{center}
\begin{tikzcd}
X \arrow[r, "u"]
& Y \arrow[r,"v"]
& Z \arrow[r, dashed, "\d"] & \text{}
\end{tikzcd} 
\end{center}
is an $\mathbb{E}$-triangle for some $\d \in \mathbb{E}(Z,X).$
\end{definition}

\begin{theorem}\label{mainTheorem} Let $(\mathscr{C}, \mathbb{E}, \mathfrak{s})$ be an extriangulated category. Then $\mathscr{C}$ has an $\mathbb{E}$-compatible triangulated structure $(\mathscr{C},T,\mathcal{T})$ if and only if for every object $ X \in \mathscr{C}$, the morphism $X \rightarrow 0$ is an $\mathbb{E}$-inflation and the morphism $0 \rightarrow X$ is an $\mathbb{E}$-deflation.
\end{theorem}
\textbf{Proof of the forward direction of Theorem \ref{mainTheorem}}
\begin{proof}
Let $(\mathscr{C}, \mathbb{E}, \mathfrak{s})$ be an extriangulated category with an $\mathbb{E}$-compatible triangulated structure $(\mathscr{C},T, \mathcal{T})$. Let $f \colon X \rightarrow Y$ be any morphism in $\mathscr{C}$. By axiom (TR1), there is a distinguished triangle, 
\begin{center}
\begin{tikzcd}
X \arrow[r, "f"]
& Y \arrow[r,"g"]
& Z \arrow[r,"h"] 
& TX \in \mathcal{T}.
\end{tikzcd} 
\end{center}
Since the triangulated structure is $\mathbb{E}$-compatible, we have that 
\begin{center}
\begin{tikzcd}
X \arrow[r, "f"]
& Y \arrow[r,"g"]
& Z \arrow[r, dashed, "\d"] & \text{}
\end{tikzcd} 
\end{center}
is an $\mathbb{E}$-triangle for some $\d \in \mathbb{E}(Z,X).$ So $f$ is an $\mathbb{E}$-inflation. In particular, we have that for every object $X \in \mathscr{C}$, the morphism $X \rightarrow 0$ is an $\mathbb{E}$-inflation.

Since $T$ is an auto-equivalence, there exists a functor $S \colon \mathscr{C} \rightarrow \mathscr{C}$ such that $S \circ T \cong 1_{\mathscr{C}}$ and $T \circ S \cong 1_{\mathscr{C}}$. Let $X$ be an arbitrary object in $\mathscr{C}$, since $T$ is essentially surjective, there exists $A \in \mathscr{C}$ such that $T(A) \cong X$. 
By axiom axiom (TR1) and (TR2), we have the following distinguished triangle
\begin{center}
\begin{tikzcd}
A \arrow[r] 
& 0 \arrow[r] 
& T(A) \arrow[r,"-1_{T(A)}"] 
&T(A)
\end{tikzcd}
\end{center}
which is isomorphic to the sextuple 
\begin{equation}\label{trivTriangle}
\begin{tikzcd}
A \arrow[r] 
& 0 \arrow[r] 
& T(A) \arrow[r,"1_{T(A)}"] 
&T(A) 
\end{tikzcd}
\end{equation}
Since $\mathcal{T}$ is closed under isomorphism by (TR1), we have that the sextuple (\ref{trivTriangle}) is also a distinguished triangle.
Since $T(A) \cong X, A \cong S(X)$ and $X \cong T \circ S (X)$ we have the following isomorphisms respectively, $\varrho \colon T(A) \rightarrow X, \pi \colon A \rightarrow S(X)$ and $\varphi \colon X \rightarrow T \circ S (X)$. Consider the following diagram. 

\begin{center}
\begin{tikzcd}
A \arrow[r] \arrow[d,"\pi"]
& 0 \arrow[r] \arrow[d,equal]
& T(A) \arrow[r,"1_{T(A)}"] \arrow[d,"\varrho"]
&T(A) \arrow[d,"\varphi \circ \varrho"]
\\ 
S(X) \arrow[r]
& 0 \arrow[r]
& X \arrow[r,"\varphi"]
& TS(X) \
\end{tikzcd}
\end{center}
Clearly the above diagram is an isomorphism of sextuples. The top row is a distinguished triangle, since $\mathcal{T}$ is closed under isomorphisms, we  have that the bottom row is also a distinguished triangle. Since $(\mathscr{C},T,\mathcal{T})$ is $\mathbb{E}$-compatible, we have that 
\begin{center}
\begin{tikzcd}
S(X) \arrow[r]
& 0 \arrow[r]
& X \arrow[r,dashed,"\d"]
& \text{}
\end{tikzcd}
\end{center}
is an $\mathbb{E}$-triangle for some $\d \in \mathbb{E}(X,S(X))$, that is to say $0 \rightarrow X$ is an $\mathbb{E}$-deflation.  
\end{proof}

The rest of this paper will be dedicated to proving the other direction. For the rest of the paper, suppose that for every object $X$ in $\mathscr{C}$, the morphism $X \rightarrow 0$ is an $\mathbb{E}$-inflation, and the morphism $0 \rightarrow X$ is an $\mathbb{E}$-deflation. We will construct an auto-equivalence $\Sigma \colon \mathscr{C} \tilde{\longrightarrow} \mathscr{C}$. Let $X$ be any object in $\mathscr{C}$. Then the morphism $X \rightarrow 0$ is an $\mathbb{E}$-inflation so there is an $\mathbb{E}$-triangle given by
\begin{equation}\label{SigmaX}
 \begin{tikzcd}
X \arrow[r]
& 0 \arrow[r]
& Z \arrow[r, dashed, "\d_{X}"] & \text{.}
\end{tikzcd}
\end{equation}
The object $Z=\text{Cone}(X \rightarrow 0)$, so it is unique up to isomorphism, hence we may choose a cone and denote it by $\Sigma X$ and fix the $\mathbb{E}$-triangle 
\begin{center}
\begin{tikzcd}
X \arrow[r]
& 0 \arrow[r]
& \Sigma X \arrow[r, dashed, "\d_{X}"] & \text{.}
\end{tikzcd}
\end{center}
Let $X$ and $Y$ be a pair of objects, and let \begin{equation}\label{extriangleX}
\begin{tikzcd}
X \arrow[r]
& 0 \arrow[r]
& \Sigma X \arrow[r, dashed, "\d_{X}"] & \text{}
\end{tikzcd}
\end{equation} and
\begin{equation}\label{extriangleY}
\begin{tikzcd}
Y \arrow[r]
& 0 \arrow[r]
& \Sigma Y\arrow[r, dashed, "\d_{Y}"] & \text{}
\end{tikzcd}
\end{equation} be the corresponding $\mathbb{E}$-triangles as above. Then by Proposition \ref{longExact} applied to the $\mathbb{E}$-triangle (\ref{extriangleX})
the following sequence is exact.
\begin{center}
\begin{tikzcd}
\mathscr{C}(\Sigma X,Y) \arrow[r]
& \mathscr{C}(0,Y) \arrow[r]
& \mathscr{C}(X,Y) \arrow[r, "(\d^{\#}_{X})_{Y}"]
& \mathbb{E}(\Sigma X,Y) \arrow[r]
& \mathbb{E}(0,Y) \arrow[r]
& \mathbb{E}(X,Y)
\end{tikzcd}
\end{center}
Since $\mathscr{C}(0,Y)=\mathbb{E}(0,Y)={0}$ in $Ab$, we have that $(\d^{\#}_{X})_{Y} \colon \mathscr{C}(X,Y) \rightarrow \mathbb{E}(\Sigma X,Y) \text{ ; } f \mapsto f_{*}\d_{X}$ is a group isomorphism. Applying Proposition \ref{longExact} to the $\mathbb{E}$-triangle (\ref{extriangleY})
the following sequence is exact in $Ab$.
\begin{center}
\begin{tikzcd}
\mathscr{C}(\Sigma X,Y) \arrow[r]
& \mathscr{C}(\Sigma X,0) \arrow[r]
& \mathscr{C}(\Sigma X,\Sigma Y) \arrow[r, "(\d_{Y\#})_{\Sigma X}"]
& \mathbb{E}(\Sigma X,Y) \arrow[r]
& \mathbb{E}(\Sigma X,0) \arrow[r]
& \mathbb{E}(\Sigma X,\Sigma Y)
\end{tikzcd}
\end{center} 
Since $\mathscr{C}(\Sigma X,0)=\mathbb{E}(\Sigma X,0)={0}$ in $Ab$, we have that $(\d_{Y\#})_{\Sigma X} \colon \mathscr{C}(\Sigma X,\Sigma Y) \rightarrow \mathbb{E}(\Sigma X,Y) \text{ ; } g \mapsto g^{*}\d_{Y}$ is a group isomorphism.

Consider the solid part of the following commutative diagram. 
\begin{equation}\label{Sigmaf}
\begin{tikzcd}
X \arrow[r] \arrow[d, "f"]
& 0 \arrow[d, equal] \arrow[r]
& \Sigma X \arrow[d, dashed, "\Sigma f"] \arrow[r,dashed,"\d_{X}"] 
& \text{}\\
 Y \arrow[r]
&0 \arrow[r]
&\Sigma Y \arrow[r,dashed,"\d_{Y}"]
& \text{}
\end{tikzcd}
\end{equation}
By the axiom (ET3),  there exists a morphism $\Sigma f \colon \Sigma X \rightarrow \Sigma Y$ such that $f_{*}\d_{X} =(\d^{\#}_{X})_{Y}(f) =(\Sigma f)^{*}\d_{Y}=(\d_{Y\#})_{\Sigma X}(\Sigma f).$ Since $(\d^{\#}_{X})_{Y}$ and $(\d_{Y\#})_{\Sigma X}$ are isomorphisms, $\Sigma f$ is the unique morphism fulfilling (ET3) in this situation. By this uniqueness we have that the $\Sigma 1_X$ is the unique morphism in $\mathscr{C}(\Sigma X, \Sigma X)$ such that $\d_{X} = (\Sigma 1_X)^{*}\d_{X}$. Therefore, since  $(1_{\Sigma X})^{*}\d_{X} = \d_{X}$, we have that $\Sigma 1_X = 1_{\Sigma X}$. Let $f \colon X \rightarrow Y$ and $g \colon Y \rightarrow Z$ be morphisms. Consider the following commutative diagram.
\begin{center}
\begin{tikzcd}
X \arrow[r] \arrow[d, "f"]
& 0 \arrow[d, equal] \arrow[r]
& \Sigma X \arrow[d, dashed, "\Sigma f"]  \arrow[r,dashed, "\d_{X}"]
& \text{} \\
 Y \arrow[r] \arrow[d,"g"]
&0 \arrow[r] \arrow[d,equal]
&\Sigma Y \arrow[d, "\Sigma g"] \arrow[r,dashed,"\d_{Y}"]
& \text{}
\\
 Z \arrow[r]
&0 \arrow[r] 
&\Sigma Z \arrow[r,dashed,"\d_{Z}"]
& \text{}
\end{tikzcd}
\end{center}
We see that $(g \circ f)_{*}\d_{X} = (\Sigma g \circ \Sigma f)^{*}\d_{Z}.$ Since by definition $\Sigma(g \circ f)$ is the unique morphism in $\mathscr{C}(\Sigma X, \Sigma Z)$ such that  $(g \circ f)_{*}\d_{X} = \Sigma(g \circ f)^{*}\d_{Z}$, we have that $ \Sigma g \circ \Sigma f = \Sigma(g \circ f)$. 

\begin{definition} Let $(\mathscr{C}, \mathbb{E}, \mathfrak{s})$ be an extriangulated category such that for every object $X \in \mathscr{C}$, the morphism $X \rightarrow 0$ is an inflation and the morphism $0 \rightarrow X$ is a deflation. Let $\Sigma \colon \mathscr{C} \rightarrow \mathscr{C}$ be the map defined as follows. For an object $X$ in $\mathscr{C}$, let $\Sigma X$ be a cone of the morphism $X \rightarrow 0$ in the $\mathbb{E}$-triangle of $\d_{X}$ as in (\ref{SigmaX}). For a morphism $f \colon X \rightarrow Y$ in $\mathscr{C}$, let $\Sigma f \colon \Sigma X \rightarrow \Sigma Y$ be the unique morphism such that $f_{*}\d_{X} = (\Sigma f)^{*}\d_{Y}$ as in (\ref{Sigmaf}).
\end{definition}

The map $\Sigma \colon \mathscr{C} \rightarrow \mathscr{C}$ is a functor by the above exposition.
% It maps an object $X$ to objects $\Sigma X$ and maps a morphism $f \colon X \rightarrow Y$ to the morphism $\Sigma f \colon \Sigma X \rightarrow \Sigma Y$, where $\Sigma X, \Sigma Y$ and $\Sigma f$ are as defined above. 
\begin{remark}
In the construction of the functor $\Sigma$, we have to make a choice about what $\Sigma X$ is since the cone is only unique up to isomorphism. However the functor $\Sigma$ is essentially independent of this choice. More precisely, the following proposition is true. 
\end{remark}
\begin{proposition} Let $(\mathscr{C}, \mathbb{E}, \mathfrak{s})$ be an extriangulated category such that for every object $X \in \mathscr{C}$, the morphism $X \rightarrow 0$ is an inflation and the morphism $0 \rightarrow X$ is a deflation. Let $\Sigma \colon \mathscr{C} \rightarrow \mathscr{C}$ and $\Sigma^{\prime} \colon \mathscr{C} \rightarrow \mathscr{C}$ be two functors defined as above, where possibly different choices of $\Sigma X$ and $\Sigma^{\prime} X$ are made for each object $X$. Then $\Sigma$ and $\Sigma^{\prime}$ are naturally isomorphic. 

\begin{proof}
Let $X,Y$ be objects in $\mathscr{C}$. By definition $\Sigma X$ and $\Sigma^{\prime} X$ are cones of $X \rightarrow 0$ and likewise $\Sigma Y$ and $\Sigma^{\prime} Y$ are cones of $Y \rightarrow 0$.
Let $f \colon X \rightarrow Y$ be a morphism in $\mathscr{C}$. Then we have the following commutative diagram of morphisms of $\mathbb{E}$-triangles,
\begin{center}
\begin{tikzcd}
X \arrow[r] \arrow[d,equal]
& 0 \arrow[r] \arrow[d, equal]
& \Sigma X \arrow[r,dashed,"\d_X "] \arrow[d, dashed, "\eta_{X}"]
& \text{}
\\
X \arrow[r] \arrow[d,"f"]
& 0 \arrow[r] \arrow[d,equal]
& \Sigma^{\prime} X \arrow[r,dashed,"\d_X "] \arrow[d,"\Sigma^{\prime} f "]
& \text{}
\\
Y \arrow[r]
& 0 \arrow[r]
&\Sigma^{\prime} Y \arrow[r,"\d_Y"]
& 
\text{}
\end{tikzcd}
\end{center}
where $\eta_{X} \colon \Sigma X \rightarrow \Sigma^{\prime} X$ is obtained by an application of the axiom (ET3) to the top two rows. By Lemma \ref{2out3iso}, the morphism $\eta_{X}$ is an isomorphism. Moreover $\d_X = \eta_{X}^{*}\d_X$ and $f_{*}\d_X = (\Sigma^{\prime} f)^{*}\d_Y$, so $f_{*}\d_X = f_{*}\eta_{X}^{*}\d_X = \eta_{X}^{*}f_{*}\d_{X} = \eta_{X}^{*}(\Sigma^{\prime}f)^{*}\d_{Y} = (\Sigma^{\prime} f \circ \eta_{X})^{*}\d_Y.$
We similarly have the following commutative diagram.
\begin{center}
\begin{tikzcd}
X \arrow[r] \arrow[d,"f"]
& 0 \arrow[r] \arrow[d, equal]
& \Sigma X \arrow[r,dashed,"\d_X "] \arrow[d, "\Sigma f"]
& \text{}
\\
Y \arrow[r] \arrow[d,equal]
& 0 \arrow[r] \arrow[d,equal]
& \Sigma^{\prime} Y \arrow[r,dashed,"\d_Y "] \arrow[d,dashed,"\eta_{Y}"]
& \text{}
\\
Y \arrow[r]
& 0 \arrow[r]
&\Sigma^{\prime} Y \arrow[r,"\d_Y"]
& 
\text{}
\end{tikzcd}
\end{center}
By Lemma \ref{2out3iso}, the morphism $\eta_{Y}$ is an isomorphism. Moreover $f_{*}\d_X = (\Sigma f)^{*}\d_Y$ and $\d_Y = (\eta_Y)^{*}\d_Y$, so $f_{*}\d_X = (\Sigma f)^{*}\eta_Y^{*}\d_Y = (\eta_{Y} \circ \Sigma f)^{*}\d_Y$. 

Since $$f_{*}\d_X = (\eta_{Y} \circ \Sigma f)^{*}\d_Y = (\Sigma^{\prime} f \circ \eta_{X})^{*}\d_Y,$$ we have that $(\eta_{Y} \circ \Sigma f) = (\Sigma^{\prime} f \circ \eta_{X})$ by the uniqueness as in (\ref{Sigmaf}). So we can define a natural transformation $\eta \colon \Sigma \Rightarrow \Sigma^{\prime}$, where for each object $X$ in $\mathscr{C}$, we associate the isomorphism $\eta_X \colon \Sigma X \rightarrow \Sigma^{\prime} X$ as obtained above. For every morphism $f \colon X \rightarrow Y$ in $\mathscr{C}$, we have that $(\eta_{Y} \circ \Sigma f) = (\Sigma^{\prime} f \circ \eta_{X})$. So $\eta$ defines a natural isomorphism, therefore $\Sigma$ and $\Sigma^{\prime}$ are naturally isomorphic.
\end{proof}
\end{proposition}

\begin{proposition} \label{SigmaAutoEquiv} Let $(\mathscr{C}, \mathbb{E}, \mathfrak{s})$ be an extriangulated category such that for every object $X \in \mathscr{C}$, the morphism $X \rightarrow 0$ is an inflation and the morphism $0 \rightarrow X$ is a deflation. Let $\Sigma \colon \mathscr{C} \rightarrow \mathscr{C}$ be as defined above. Then $\Sigma \colon \mathscr{C} \rightarrow \mathscr{C}$ is an additive auto-equivalence. 

\begin{proof}
We first show that $\Sigma$ is an additive functor. We can see from the $\mathbb{E}$-triangle 
\begin{center}
\begin{tikzcd}
0 \arrow[r]
& 0 \arrow[r]
& \Sigma 0 \arrow[r, dashed, "\d_{0}"] & \text{,}
\end{tikzcd}
\end{center}
that $\Sigma 0 \cong 0$ since Cone$(0 \rightarrow 0)=0$, and is unique up to isomorphism.

 Let $X,Y$ be any pair of objects and consider the $\mathbb{E}$-triangle 
\begin{center} 
 \begin{tikzcd}
X \oplus Y \arrow[r]
& 0 \arrow[r]
& \Sigma (X \oplus Y) \arrow[r, dashed, "\d_{X \oplus Y}"] & \text{.}
\end{tikzcd}
\end{center}
Since $\mathfrak{s}$ is an additive realisation, we have that $\mathfrak{s}(\d_X \oplus \d_Y) = \mathfrak{s}(\d_X) \oplus \mathfrak{s}(\d_Y)$. So the following is an $\mathbb{E}$-triangle $$\begin{tikzcd}
X \oplus Y \arrow[r]
& 0 \arrow[r]
& \Sigma X \oplus \Sigma Y \arrow[r, dashed, "\d_{X} \oplus \d_{Y}"] & \text{.} \end{tikzcd}$$
Observe that the solid part of the following diagram commutes,
\begin{center}
\begin{tikzcd}
X \oplus Y \arrow[r] \arrow[d, equal]
& 0 \arrow[d, equal] \arrow[r]
& \Sigma(X \oplus Y) \arrow[d, dashed, "c"] \arrow[r,dashed,"\d_{X \oplus Y}"]
& \text{} \\
  X \oplus Y\arrow[r]
&0 \arrow[r]
&\Sigma X \oplus \Sigma Y \arrow[r,dashed,"\d_{X} \oplus \d_{Y}"]
& \text{} \end{tikzcd}
\end{center}
so by the axiom (ET3), there exists a morphism $c \colon \Sigma (X \oplus Y) \rightarrow \Sigma X \oplus \Sigma Y$, such that the above diagram is a morphism of $\mathbb{E}$-triangles. By Lemma \ref{2out3iso}, $c$ is an isomorphism, in other words $\Sigma (X \oplus Y) \cong \Sigma X \oplus \Sigma Y$. Since $\Sigma$ is a functor between additive categories, this shows that $\Sigma$ is an additive functor. 

All that is left is to show $\Sigma$ is an equivalence. To this end, we will show that $\Sigma$ is essentially surjective and fully faithful. We start by showing it is essentially surjective. Let $Y$ be any object in $\mathscr{C}$. By assumption the morphism $0 \rightarrow Y$ is a deflation, so there is an $\mathbb{E}$-triangle 
\begin{center}
\begin{tikzcd}
X \arrow[r]
& 0 \arrow[r]
& Y \arrow[r, dashed, "\d"] & \text{.} \end{tikzcd}
\end{center}
Observe that the solid part of the following diagram commutes.
\begin{center}
\begin{tikzcd}
X \arrow[r] \arrow[d, equal]
& 0 \arrow[d, equal] \arrow[r]
& \Sigma X \arrow[d, dashed, "c"] \arrow[r,dashed,"\d_{X}"]
& \text{} \\
  X \arrow[r]
&0 \arrow[r]
& Y \arrow[r,dashed,"\d"]
& \text{} \end{tikzcd}
\end{center}
So by (ET3) and Lemma \ref{2out3iso}, there exists an isomorphism $c \colon \Sigma X \rightarrow Y$, in other words $Y \cong \Sigma X$, so $\Sigma$ is essentially surjective. 

Let $X,Y$ be a pair of objects in $\mathscr{C}$.  Consider the map $\Sigma_{X,Y} \colon \mathscr{C}(X,Y) \rightarrow \mathscr{C}(\Sigma X, \Sigma Y)$ where $\Sigma_{X,Y}(f) = \Sigma f$. Suppose $\Sigma f = \Sigma f^{\prime}$, then from (\ref{Sigmaf}) we have that $$f_{*}\d_{X} = (\Sigma f)^{*}d_{Y} = (\Sigma f^{\prime})^{*}\d_{Y} = f^{\prime}_{*}\d_{X}.$$ so $$(\d_{X}^{\#})_{Y}(f) = (\d_{X}^{\#})_{Y}(f^{\prime}).$$ Since $(\d_{X}^{\#})_{Y}$ is an isomorphism, we have that $f = f^{\prime}$, that is to say $\Sigma_{X,Y}$ is injective. 

Let $g \colon \Sigma X \rightarrow \Sigma Y$ be any morphism in $\mathscr{C}(\Sigma X, \Sigma Y)$. From the $\mathbb{E}$-triangle 
\begin{center}
\begin{tikzcd}
Y \arrow[r]
& 0 \arrow[r]
& \Sigma Y \arrow[r, dashed, "\d_{Y}"] & \text{,} \end{tikzcd}
\end{center}
we obtain the $\mathbb{E}$-triangle 
\begin{center}
\begin{tikzcd}
Y \arrow[r]
& M \arrow[r]
& \Sigma X \arrow[r, dashed, "g^{*}\d_{Y}"] & \text{.} \end{tikzcd}
\end{center}
We can then construct the solid part of the following commutative diagram, where the triple of morphisms $(1_Y, 0, g)$ between the second and third row is from an application of Lemma \ref{*operation}. 
\begin{center}
\begin{tikzcd}
X \arrow[r] \arrow[d,dashed, "f"]
& 0 \arrow[r] \arrow[d]
& \Sigma X \arrow[r, dashed, "\d_{X}"] \arrow[d,equal]
& \text{} \\
Y \arrow[r] \arrow[d, equal]
&M \arrow[r] \arrow[d]
&\Sigma X \arrow[r, dashed, "g^{*}\d_{Y}"] \arrow[d,"g"]
& \text{}  \\
Y \arrow[r]
& 0 \arrow[r]
& \Sigma Y \arrow[r, dashed, "d_{Y}"]
& \text{}
\end{tikzcd}
\end{center}
By applying the axiom $\text{(ET3)}^{\text{op}}$ to the solid commutative diagram given by the solid square between the first and second row; we have that there exists a morphism $f \colon X \rightarrow Y$ such that the top rectangle commutes and $f_{*}\d_{X} = g^{*}\d_{Y}$. Since $\Sigma f$ is the unique morphism in $\mathscr{C}(\Sigma X, \Sigma Y)$ such that $f_{*}\d_{X} = (\Sigma f)^{*}\d_{Y}$, we have that $g = \Sigma f$. Therefore the map $\Sigma_{X,Y}$ is surjective as well as injective. So we conclude that $\Sigma$ is fully faithful. This completes the proof.
\end{proof}
\end{proposition}

Now let $\mathbf{E}^{1} \colon \mathscr{C}^{\text{op}} \times \mathscr{C} \rightarrow Ab$ be the bifunctor defined by $\mathbf{E}^{1}(-,-) := \mathscr{C}(-,\Sigma -).$ We will show that it is a biadditive functor. 

\begin{lemma}\label{lemmaE1} Let $(\mathscr{C}, \mathbb{E}, \mathfrak{s})$ be an extriangulated category such that for every object $A \in \mathscr{C}$, the morphism $A \rightarrow 0$ is an inflation and the morphism $0 \rightarrow A$ is a deflation. Let $X,Y$ be a pair of objects in $\mathscr{C}$. Then $\mathbf{E}^{1}(X,Y):=\mathscr{C}(X,\Sigma Y) \cong \mathbb{E}(X,Y)$. 
\begin{proof} By Proposition \ref{longExact} applied to the $\mathbb{E}$-triangle \begin{tikzcd}
Y \arrow[r]
& 0 \arrow[r]
& \Sigma Y \arrow[r, dashed, "\d_{Y}"] & \text{,} \end{tikzcd}
the following sequence is exact. 
\begin{center}
\begin{tikzcd}
\mathscr{C}(X,Y) \arrow[r]
& \mathscr{C}(X,0) \arrow[r]
& \mathscr{C}(X,\Sigma Y) \arrow[r, "(\d_{Y\#})_{X}"]
& \mathbb{E}(X,Y) \arrow[r]
& \mathbb{E}(X,0) \arrow[r]
& \mathbb{E}(X,\Sigma Y)
\end{tikzcd}
\end{center} 
Since $\mathscr{C}(X,0)=\mathbb{E}(X,0)={0}$ we have $(\d_{Y\#})_{X} \colon \mathscr{C}(X,\Sigma Y) \rightarrow \mathbb{E}(X,Y)$ is an isomorphism.  
\end{proof}
\end{lemma}

\begin{lemma}\label{ET1prf} The functor $\mathbf{E}^{1} \colon \mathscr{C}^{\text{op}} \times \mathscr{C} \rightarrow Ab$ is a biadditive functor. 
\begin{proof}
From Lemma \ref{lemmaE1}, we have that $\mathbf{E}^{1}(X,0) \cong \mathbb{E}(X,0) \cong {0} \cong \mathbb{E}(0,X) \cong \mathbf{E}^{1}(0,X)$ and $\mathbf{E}^{1}(X,Y \oplus Z) \cong \mathbb{E}(X,Y \oplus Z) \cong \mathbb{E}(X,Y) \oplus \mathbb{E}(X,Z) \cong \mathbf{E}^{1}(X,Y) \oplus \mathbf{E}^{1}(X,Z)$. This shows $\mathbf{E}^{1}$ is additive in the second argument. Dually it is also additive in the first argument, therefore $\mathbf{E}^{1}$ is a biadditive functor.
\end{proof}
\end{lemma}
 
\begin{lemma} Let $(\mathscr{C}, \mathbb{E}, \mathfrak{s})$ be an extriangulated category such that for every object $C \in \mathscr{C}$, the morphism $C \rightarrow 0$ is an inflation and the morphism $0 \rightarrow C$ is a deflation. Let $A$ be an object in $\mathscr{C}$ and $f \colon X \rightarrow Y$ be a morphism in $\mathscr{C}$. Then we have that $\mathbf{E}^{1}(A,f) \cong \mathbb{E}(A,f)$ and $\mathbf{E}(f^{\text{op}},A) \cong \mathbb{E}(f^{\text{op}},A)$ in the category Mor$(Ab)$, the category of morphisms of $Ab$.
\begin{proof}
By Lemma \ref{lemmaE1}, we have that the map $(\d_{X\#})_{A} : \mathbf{E}^{1}(A,X) \rightarrow \mathbb{E}(A,X)$; where $(\d_{X\#})_{A}(\varepsilon) = \varepsilon^{*}\d_{X},$ and the map $(\d_{Y\#})_{A} : \mathbf{E}^{1}(A,Y) \rightarrow \mathbb{E}(A,Y)$; where $(\d_{Y\#})_{A}(\varepsilon) = \varepsilon^{*}\d_{Y}$, are group isomorphisms. So we just need to show that the following diagram commutes. 

\begin{center}
\begin{tikzcd}[column sep=3.5em,row sep=3.5em]
\mathbf{E}^{1}(A,X) \arrow[r, "\mathbf{E}^{1}(A{,}f)"] \arrow[d,"(\d_{X\#})_{A}" left]
& \mathbf{E}^{1}(A,Y) \arrow[d,"(\d_{Y\#})_{A}"] \\
\mathbb{E}(A,X) \arrow[r, "\mathbb{E}(A{,}f)"]
&\mathbb{E}(A,Y)
\end{tikzcd}
\end{center} 
Take any $\varepsilon \colon A \rightarrow \Sigma X$ in $\mathbf{E}^{1}(A,X)$. Then $(\d_{X\#})_{A}(\varepsilon) = \varepsilon^{*}\d_{X}$, therefore $\mathbb{E}(A,f)\circ(\d_{X\#})_{A}(\varepsilon) = f_{*}\varepsilon^{*}\d_{X}.$ On the other hand $\mathbf{E}^{1}(f,A)(\varepsilon) = \Sigma f \circ \varepsilon$, therefore $(\d_{Y\#})_{A} \circ \mathbf{E}^{1}(A,f)(\varepsilon) = (\Sigma f \circ \varepsilon)^{*}\d_{Y} = \varepsilon^{*} (\Sigma f)^{*} \d_{Y}$. Since $f_{*}\d_{X} = (\Sigma f)^{*}\d_{Y}$ we have that $\varepsilon^{*}(\Sigma f)^{*}\d_{Y} = \varepsilon^{*}f_{*} \d_{X} = f_{*}\varepsilon^{*} \d_{X}$. So the diagram commutes, that is $\mathbf{E}^{1}(A,f) \cong \mathbb{E}(A,f)$ in Mor$(Ab)$. The argument for the other statement is dual. 
\end{proof}
\end{lemma}

We now define a correspondence $\mathfrak{r}$ for the category $\mathscr{C}$ endowed with the biadditive functor $\mathbf{E}^{1}$, which will associate an equivalence class $\mathfrak{r}(\varepsilon)$ to any extension $\varepsilon \in \mathbf{E}^{1}(C,A)$. Let $A,C$ be a pair of objects in $\mathscr{C}$, then by the isomorphism in Lemma \ref{lemmaE1}, $\mathbf{E}^{1}(C,A) \cong \mathbb{E}(C,A)$, so every $\varepsilon$ corresponds to $\varepsilon^{*}\d_{A} \in \mathbb{E}(C,A)$. So we set $\mathfrak{r}(\varepsilon) := \mathfrak{s}(\varepsilon^{*}\d_{A})$. 

\begin{proposition}\label{ET2prf} Let $(\mathscr{C}, \mathbb{E}, \mathfrak{s})$ be an extriangulated category such that for every object $X \in \mathscr{C}$, the morphism $X \rightarrow 0$ is an inflation and the morphism $0 \rightarrow X$ is a deflation. Let $\mathfrak{r}$ be the correspondence which associates the equivalence class $\mathfrak{r}(\varepsilon) = \mathfrak{s}(\varepsilon^{*}\d_{A})$ to any $\mathbf{E}^{1}$-extension $\varepsilon \in \mathbf{E}^{1}(C,A)$, for any pair of objects $A,C$ in $\mathscr{C}$. Then $\mathfrak{r}$ is an additive realisation of $\mathbf{E}^{1}$. 

\begin{proof} Let $\varepsilon \in \mathbf{E}^{1}(C,A)$ and $\varepsilon^{\prime} \in \mathbf{E}^{1}(C^{\prime},A^{\prime})$ be $\mathbf{E}^{1}$-extensions with $\mathfrak{r}(\varepsilon) = [A \overset{x}{\longrightarrow} B \overset{y}{\longrightarrow} C]$ and $\mathfrak{r}(\varepsilon^{\prime}) = [A^{\prime} \overset{x^{\prime}}{\longrightarrow} B^{\prime} \overset{y^{\prime}}{\longrightarrow} C^{\prime}]$. Suppose we have a morphism of $\mathbf{E}^{1}$-extensions $(a,c) \colon \varepsilon \rightarrow \varepsilon^{\prime}$, that is to say $\mathbf{E}^{1}(C,a)(\varepsilon)=\Sigma a \circ \varepsilon = \mathbf{E}^{1}(c,A^{\prime})(\varepsilon^{\prime})= \varepsilon^{\prime} \circ c$.  Since $\mathfrak{r}(\varepsilon) = \mathfrak{s}(\varepsilon^{*}\d_{A})$ and $\mathfrak{r}(\varepsilon^{\prime}) = \mathfrak{s}((\varepsilon^{\prime})^{*}\d_{A^{\prime}})$, we have the following diagram of $\mathbb{E}$-triangles. 
\begin{center}
\begin{tikzcd}
A \arrow[r, "x"] \arrow[d,"a"]
& B \arrow[r,"y"] \arrow[d, dashed, "b"]
& C \arrow[r, dashed, "\varepsilon^{*}\d_{A}"] \arrow[d,"c"] & \text{} \\
A^{\prime} \arrow[r, "x^{\prime}"]
& B^{\prime} \arrow[r,"y^{\prime}"]
& C^{\prime} \arrow[r, dashed, "(\varepsilon^{\prime})^{*}\d_{A^{\prime}}"] &\text{}
\end{tikzcd}
\end{center}
Recall from (\ref{Sigmaf}) that $a_{*}\d_{A} = (\Sigma a)^{*}\d_{A^{\prime}}$. Therefore $a_{*}\varepsilon^{*}\d_{A} = \varepsilon^{*}a_{*}\d_{A} = \varepsilon^{*}(\Sigma a)^{*}\d_{A^{\prime}} = (\Sigma a \circ \varepsilon)^{*}\d_{A^{\prime}} = (\varepsilon^{\prime} \circ c)^{*}\d_{A^{\prime}}= c^{*} (\varepsilon^{\prime})^{*}\d_{A^{\prime}}$, in other words $(a,c) \colon \varepsilon^{*}\d_{A} \rightarrow (\varepsilon^{\prime})^{*}\d_{A^{\prime}}$ is a morphism of $\mathbb{E}$-extensions. Therefore, since $\mathfrak{s}$ is a realisation, there exists a morphism $b \colon B \rightarrow B^{\prime}$ making the above diagram commute. This shows that $\mathfrak{r}$ is a realisation. 

What is left is to show that $\mathfrak{r}$ is an additive realisation. For any objects $A,C$ in $\mathscr{C}$, take split $\mathbf{E}^{1}$-extension $ 0 \in \mathbf{E}^{1}(C,A)$, then $$\mathfrak{r}(0)=\mathfrak{s}(0^{*}\d_{A}) = \mathfrak{s}(0)=0.$$
Let $\varepsilon \in \mathbf{E}^{1}(C,A)$ and $\varepsilon^{\prime} \in \mathbf{E}^{1}(C^{\prime},A^{\prime})$ be $\mathbf{E}^{1}$-extensions with $\mathfrak{r}(\varepsilon) = \mathfrak{s}(\varepsilon^{*}\d_{A}) = [A \overset{x}{\longrightarrow} B \overset{y}{\longrightarrow} C]$ and $\mathfrak{r}(\varepsilon^{\prime}) = \mathfrak{s}(\varepsilon^{\prime *}\d_{A^{\prime}}) =[ A^{\prime} \overset{x^{\prime}}{\longrightarrow} B^{\prime} \overset{y^{\prime}}{\longrightarrow} C^{\prime}]$. By Definition \ref{additiveRealisation}, we need to show that $$\mathbf{E}^{1}(p_C,i_A)(\varepsilon) + \mathbf{E}^{1}(p_{C^{\prime}},i_{A^{\prime}})(\varepsilon^{\prime})$$ is realised by the direct sum.
\begin{center}
\begin{tikzcd}
A \oplus A^{\prime} \arrow[r,"x \oplus x^{\prime}"]
& B \oplus B^{\prime} \arrow[r,"y \oplus y^{\prime}"]
& C \oplus C^{\prime}.
\end{tikzcd}
\end{center}
First observe that by the definition of $ \Sigma$ on morphisms, $(i_{A})_{*}\d_{A} = (\Sigma i_{A})^{*}\d_{A \oplus A^{\prime}}$ due to the following commutative diagram.
\begin{equation}
\begin{tikzcd}
A \arrow[r] \arrow[d, "i_{A}"]
& 0 \arrow[d, equal] \arrow[r]
& \Sigma A \arrow[d, dashed, "\Sigma i_{A}"] \arrow[r,dashed,"\d_{A}"] 
& \text{}\\
 A \oplus A^{\prime} \arrow[r]
&0 \arrow[r]
&\Sigma (A \oplus A^{\prime}) \arrow[r,dashed,"\d_{A \oplus A^{\prime}}"]
& \text{}
\end{tikzcd}
\end{equation}
Likewise, we have that $(i_{A^{\prime}})_{*}\d_{A^{\prime}} = (\Sigma i_{A^{\prime}})^{*}\d_{A \oplus A^{\prime}}$.
By direct calculation, we have the following. 
\begin{align*}
(i_{A})_{*}(p_{C})^{*}(\varepsilon^{*}\d_{A}) + (i_{A^{\prime}})_{*}(p_{C^{\prime}})^{*}(\varepsilon^{\prime *}\d_{A})  &=  (p_{C})^{*}\varepsilon^{*}((i_{A})_{*}\d_{A}) +  (p_{C^{\prime}})^{*}\varepsilon^{\prime *}((i_{A^{\prime}})_{*}\d_{A}) \\
&= (\varepsilon \circ p_{C})^{*}((i_{A})_{*}\d_{A}) + (\varepsilon^{\prime} \circ p_{C^{\prime}})^{*}((i_{A^{\prime}})_{*}\d_{A^{\prime}})  \\
&= (\varepsilon \circ p_{C})^{*}((\Sigma i_{A^{\prime}})^{*}\d_{A \oplus A}) + (\varepsilon^{\prime} \circ p_{C^{\prime}})^{*}((\Sigma i_{A^{\prime}})^{*}\d_{A \oplus A^{\prime}})    \\
&=  ( \Sigma i_{A} \circ \varepsilon \circ p_{C})^{*}\d_{A \oplus A^{\prime}} + ( \Sigma i_{A^{\prime}} \circ \varepsilon^{\prime} \circ p_{C^{\prime}})^{*}\d_{A \oplus A^{\prime}}  \\
&=  ( \Sigma i_{A} \circ \varepsilon \circ p_{C} + \Sigma i_{A^{\prime}} \circ \varepsilon^{\prime} \circ p_{C^{\prime}})^{*}\d_{A \oplus A^{\prime}} \\
&= (\mathbf{E}^{1}(p_C,i_A)(\varepsilon) + \mathbf{E}^{1}(p_{C^{\prime}},i_{A^{\prime}})(\varepsilon^{\prime}))^{*}\d_{A \oplus A^{\prime}}
\end{align*}
By definition $$\mathfrak{r}(\mathbf{E}^{1}(p_C,i_A)(\varepsilon) + \mathbf{E}^{1}(p_{C^{\prime}},i_{A^{\prime}})(\varepsilon^{\prime})) = \mathfrak{s}((\mathbf{E}^{1}(p_C,i_A)(\varepsilon) + \mathbf{E}^{1}(p_{C^{\prime}},i_{A^{\prime}})(\varepsilon^{\prime}))^{*}\d_{A \oplus A^{\prime}}),$$
by the above calculation 
$$\mathfrak{r}(\mathbf{E}^{1}(p_C,i_A)(\varepsilon) + \mathbf{E}^{1}(p_{C^{\prime}},i_{A^{\prime}})(\varepsilon^{\prime})) =\mathfrak{s}((i_{A})_{*}(p_{C})^{*}(\varepsilon^{*}\d_{A}) + (i_{A^{\prime}})_{*}(p_{C^{\prime}})^{*}(\varepsilon^{\prime *}\d_{A})).$$
Since $\mathfrak{s}$ is an additive realisation, we have that $(i_{A})_{*}(p_{C})^{*}(\varepsilon^{*}\d_{A}) + (i_{A^{\prime}})_{*}(p_{C^{\prime}})^{*}(\varepsilon^{\prime *}\d_{A})$ is realised by the following direct sum with respect to $\mathfrak{s}$,
\begin{center}
\begin{tikzcd}
A \oplus A^{\prime} \arrow[r,"x \oplus x^{\prime}"]
& B \oplus B^{\prime} \arrow[r,"y \oplus y^{\prime}"]
& C \oplus C^{\prime}.
\end{tikzcd}
\end{center}
hence $$\mathfrak{r}(\mathbf{E}^{1}(p_C,i_A)(\varepsilon) + \mathbf{E}^{1}(p_{C^{\prime}},i_{A^{\prime}})(\varepsilon^{\prime})) = [A \oplus A^{\prime} \overset{x \oplus x^{\prime}}{\longrightarrow} B \oplus B^{\prime} \overset{y \oplus y^{\prime}}{\longrightarrow} C \oplus C^{\prime}].$$ This completes the proof. So $\mathfrak{r}$ is an additive realisation.
\end{proof}
\end{proposition} 

\begin{proposition}\label{ET3prf} The triple $(\mathscr{C},\mathbf{E}^{1},\mathfrak{r})$ satisfies the axioms (ET3) and $\text{(ET3)}^{\text{op}}$. 
\begin{proof} Let $\varepsilon \in \mathbf{E}^{1}(C,A)$ and $\varepsilon^{\prime} \in \mathbf{E}^{1}(C^{\prime},A^{\prime})$ be any pair of $\mathbf{E}^{1}$-extensions realised by the sequences $A \overset{x}{\longrightarrow} B \overset{y}{\longrightarrow} C$ and $A^{\prime} \overset{x^{\prime}}{\longrightarrow} B^{\prime} \overset{y^{\prime}}{\longrightarrow} C^{\prime}$ respectively. Consider the following commutative diagram in $(\mathscr{C},\mathbf{E}^{1},\mathfrak{r})$. 
\begin{center}
\begin{tikzcd}
A \arrow[r, "x"] \arrow[d, "a"]
& B \arrow[d, "b"] \arrow[r,"y"]
& C \\
 A^{\prime} \arrow[r,"x^{\prime}"]
&B^{\prime} \arrow[r,"y^{\prime}"]
&C^{\prime} \end{tikzcd}
\end{center}
Since $\mathfrak{r}(\varepsilon) = \mathfrak{s}(\varepsilon^{*}\d_{A})$ and $\mathfrak{r}(\varepsilon) = \mathfrak{s}((\varepsilon^{\prime})^{*}\d_{A^{\prime}})$, the sequences $A \overset{x}{\longrightarrow} B \overset{y}{\longrightarrow} C$ and $A^{\prime} \overset{x^{\prime}}{\longrightarrow} B^{\prime} \overset{y^{\prime}}{\longrightarrow} C^{\prime}$ realise $\varepsilon^{*}\d_{A}$ and $(\varepsilon^{\prime})^{*}\d_{A^{\prime}}$  respectively in $(\mathscr{C},\mathbb{E},\mathfrak{s})$. We also have the following diagram with the solid part commuting in $(\mathscr{C},\mathbb{E},\mathfrak{s})$. 
\begin{center}
\begin{tikzcd}
A \arrow[r, "x"] \arrow[d, "a"]
& B \arrow[d, "b"] \arrow[r,"y"]
& C \arrow[d,dashed,"c"] \\
 A^{\prime} \arrow[r,"x^{\prime}"]
&B^{\prime} \arrow[r,"y^{\prime}"]
&C^{\prime} \end{tikzcd}
\end{center}
So by (ET3) in $(\mathscr{C},\mathbb{E},\mathfrak{s})$, there exists a morphism $c \colon C \rightarrow C^{\prime}$ making the above diagram commute, such that $a_{*}(\varepsilon^{*}\d_{A}) = c^{*} (\varepsilon^{\prime})^{*}\d_{A^{\prime}}$. Recall from (\ref{Sigmaf}) that $a_{*}\d_{A} = (\Sigma a)^{*}\d_{A^{\prime}}$. Therefore $$a_{*}\varepsilon^{*}\d_{A} = \varepsilon^{*}a_{*}\d_{A} = \varepsilon^{*}(\Sigma a)^{*}\d_{A^{\prime}} = (\Sigma a \circ \varepsilon)^{*}\d_{A^{\prime}} = (\varepsilon^{\prime} \circ c)^{*}\d_{A^{\prime}}= c^{*} (\varepsilon^{\prime})^{*}\d_{A^{\prime}}.$$ Since $$(\d_{A^{\prime}\#})_{C}(\Sigma a \circ \varepsilon) = (\Sigma a \circ \varepsilon)^{*}\d_{A^{\prime}} =  (\varepsilon^{\prime} \circ c)^{*}\d_{A^{\prime}} = (\d_{A^{\prime}\#})_{C}(\varepsilon^{\prime} \circ c),$$
and $(\d_{A^{\prime}\#})_{C}$ is an isomorphism, we have that $\Sigma a \circ \varepsilon = \varepsilon^{\prime} \circ c$. Therefore, in $(\mathscr{C},\mathbf{E}^{1},\mathfrak{r})$ we have that there exists a morphism $c \colon C \rightarrow C^{\prime}$ such that the following diagram commutes, 
\begin{center}
\begin{tikzcd}
A \arrow[r, "x"] \arrow[d, "a"]
& B \arrow[d, "b"] \arrow[r,"y"]
& C \arrow[d,"c"] \\
 A^{\prime} \arrow[r,"x^{\prime}"]
&B^{\prime} \arrow[r,"y^{\prime}"]
&C^{\prime} \end{tikzcd}
\end{center}
and $\mathbf{E}^{1}(C,a)(\varepsilon) = \mathbf{E}^{1}(c,A)(\varepsilon^{\prime})$, hence $(\mathscr{C},\mathbf{E}^{1},\mathfrak{r})$ satisfies (ET3). The proof that $(\mathscr{C},\mathbf{E}^{1},\mathfrak{r})$ satisfies (ET3)$^{\text{op}}$ is dual.
\end{proof}
\end{proposition}

\begin{proposition}\label{ET4prf}  The triple $(\mathscr{C},\mathbf{E}^{1},\mathfrak{r})$ satisfies (ET4) and (ET4)$^{\text{op}}$.
\begin{proof}
 Let $\varepsilon \in \mathbf{E}^{1}(D,A)$ and $\varepsilon^{\prime} \in \mathbf{E}^{1}(F,B)$ be any pair of $\mathbf{E}^{1}$-extensions realised by the sequences, $A \overset{f}{\longrightarrow} B \overset{f^{\prime}}{\longrightarrow} D$ and $B\overset{g}{\longrightarrow} C \overset{g^{\prime}}{\longrightarrow} F$ respectively. 
Since $\mathfrak{r}(\varepsilon) = \mathfrak{s}(\varepsilon^{*}\d_{A})$ and $\mathfrak{r}(\varepsilon^{\prime}) = \mathfrak{s}((\varepsilon^{\prime})^{*}\d_{B})$, we have that the $\mathbb{E}$-extensions $\varepsilon^{*}\d_{A}$ and $(\varepsilon^{\prime})^{*}\d_{B}$ are realised by the sequences $A \overset{f}{\longrightarrow} B \overset{f^{\prime}}{\longrightarrow} D$ and $B\overset{g}{\longrightarrow} C \overset{g^{\prime}}{\longrightarrow} F$ respectively. Therefore by (ET4) applied to the $\mathbb{E}$-extensions, $\varepsilon^{*}\d_{A}$ and $(\varepsilon^{\prime})^{*}\d_{B}$, we have the following diagram of $\mathbb{E}$-triangles
\begin{center}
\begin{tikzcd}
A \arrow[r, "f"] \arrow[d, equal]
& B \arrow[d, "g"] \arrow[r,"f^{\prime}"]
& D \arrow[d,"d"] \arrow[r,dashed,"\varepsilon^{*}\d_{A}"]
& \text{} 
\\
 A \arrow[r,"h"]
&C \arrow[r,"h^{\prime}"] \arrow[d,"g^{\prime}"]
&E \arrow[d,"e"] \arrow[r,dashed,"\d^{\prime \prime}"]
& \text{}
 \\
 & F \arrow[r,equal] \arrow[d,dashed,"(\varepsilon^{\prime})^{*}\d_{B}"] & F \arrow[d,dashed,"(f^{\prime})_{*}(\varepsilon^{\prime})^{*}\d_{B}"]
 \\ & \text{}  & \text{}
\end{tikzcd}
\end{center}
in $(\mathscr{C},\mathbb{E},\mathfrak{s})$ and an $\mathbb{E}$-extension $\d^{\prime \prime} \in \mathbb{E}(E,A)$ realised by the sequence $A \overset{h}{\longrightarrow} C \overset{h^{\prime}}{\longrightarrow} E$, such that the following compatibilities are satisfied;
\begin{enumerate}[(i)]
\item $\mathfrak{s}((f^{\prime})_{*}(\varepsilon^{\prime})^{*}\d_{B}) = [D \overset{d}{\longrightarrow} E \overset{e}{\longrightarrow} F].$

\item $d^{*} \d^{\prime \prime} = \varepsilon^{*}\d_{A}.$

\item $f_{*}\d^{\prime \prime} = e^{*}(\varepsilon^{\prime})^{*}\d_{B}.$
\end{enumerate}
Since $\mathbb{E}(E,A) \cong \mathbf{E}^{1}(E,A)$ we have that there exists $\varepsilon^{\prime \prime} \in \mathscr{C}(E,\Sigma A)$ such that $\d^{\prime \prime} = (\varepsilon^{\prime \prime})^{*}\d_{A}$. Moreover since $\mathfrak{r}(\varepsilon^{\prime \prime}) = \mathfrak{s}((\varepsilon^{\prime \prime})^{*}\d_{A})$, as an $\mathbf{E}^{1}$-extension $\varepsilon^{\prime \prime}$ is realised by the sequence, $$A \overset{h}{\longrightarrow} C \overset{h^{\prime}}{\longrightarrow} E.$$

Consider the $\mathbb{E}$-extension $(f^{\prime})_{*}(\varepsilon^{\prime})^{*}\d_{B}$. Recall that $(f^{\prime})_{*}\d_{B} = (\Sigma f^{\prime})^{*}\d_{D}$ as in (\ref{Sigmaf}) therefore
\begin{equation}\label{ET4Eqn1}
(f^{\prime})_{*}(\varepsilon^{\prime})^{*}\d_{B} = (\varepsilon^{\prime})^{*}(f^{\prime})_{*}\d_{B} =  (\varepsilon^{\prime})^{*}(\Sigma f^{\prime})^{*}\d_{D} = (\Sigma f^{\prime} \circ \varepsilon^{\prime})^{*}\d_{D}.
\end{equation}
So we have by compatibility (i) and (\ref{ET4Eqn1}) that $$\mathfrak{r}(\mathbf{E}^{1}(F,f^{\prime})(\varepsilon^{\prime})) = \mathfrak{r}(\Sigma f^{\prime} \circ \varepsilon^{\prime}) = \mathfrak{s}((\Sigma f^{\prime} \circ \varepsilon^{\prime})^{*}\d_{D}) = \mathfrak{s}((f^{\prime})_{*}(\varepsilon^{\prime})^{*}\d_{B}) = [D \overset{d}{\longrightarrow} E \overset{e}{\longrightarrow} F].$$

By compatibility (ii), we have that 
$$(\d_{A\#})_{D}(\varepsilon^{\prime \prime} \circ d) = d^{*}(\varepsilon^{\prime \prime})^{*}\d_{A} =  d^{*} \d^{\prime \prime} = \varepsilon^{*}\d_{A} = (\d_{A\#})_{D}(\varepsilon).$$
Since $(\d_{A\#})_{D}$ is an isomorphism, we have that $\varepsilon^{\prime \prime} \circ d = \varepsilon$, in particular, we have that $$\mathbf{E}^{1}(d,A)(\varepsilon^{\prime \prime}) = \varepsilon.$$

Recall that $f_{*}\d_{A}  =\Sigma f ^{*}\d_{B}$ as in (\ref{Sigmaf}), using this equality and compatibility (iii) we have that 
\begin{equation}\label{ET4Eqn2} f_{*}\d^{\prime \prime} = f_{*}(\varepsilon^{\prime \prime})^{*}\d_{A} = (\varepsilon^{\prime \prime})^{*}f_{*}\d_{A} =  (\varepsilon^{\prime \prime})^{*} \Sigma f^{*}\d_{B} = (\Sigma f \circ \varepsilon^{\prime \prime})^{*}\d_{B} = (\varepsilon^{\prime} \circ e)^{*}\d_{B} = e^{*}(\varepsilon^{\prime})^{*}\d_{B}.
\end{equation}
By (\ref{ET4Eqn2}) we have that $(\Sigma f \circ \varepsilon^{\prime \prime})^{*}\d_{B} = (\varepsilon^{\prime} \circ e)^{*}\d_{B}$, equivalently $(\d_{B\#})_{E}(\Sigma f \circ \varepsilon^{\prime \prime}) = (\d_{B\#})_{E}(\varepsilon^{\prime} \circ e)$. Since $(\d_{B\#})_{E}$ is an isomorphism, we have that $\Sigma f \circ \varepsilon^{\prime \prime} = \varepsilon^{\prime} \circ e$. In particular, we have that 
$$ \mathbf{E}^{1}(E,f)(\varepsilon^{\prime \prime}) = \mathbf{E}^{1}(e,B)(\varepsilon^{\prime}).$$
To conclude, we have shown that given any pair of $\mathbf{E}^{1}$-extensions $\varepsilon \in \mathbf{E}^{1}(D,A)$ and $\varepsilon^{\prime} \in \mathbf{E}^{1}(F,B)$ realised by the sequences, $A \overset{f}{\longrightarrow} B \overset{f^{\prime}}{\longrightarrow} D$ and $B\overset{g}{\longrightarrow} C \overset{g^{\prime}}{\longrightarrow} F$ respectively. There exists an object $E$ in $(\mathscr{C},\mathbf{E}^{1},\mathfrak{r})$, a commutative diagram
\begin{center}
\begin{tikzcd}
A \arrow[r, "f"] \arrow[d, equal]
& B \arrow[d, "g"] \arrow[r,"f^{\prime}"]
& D \arrow[d,"d"] \arrow[r,dashed,"\varepsilon"]
& \text{} 
\\
 A \arrow[r,"h"]
&C \arrow[r,"h^{\prime}"] \arrow[d,"g^{\prime}"]
&E \arrow[d,"e"] \arrow[r,dashed,"\varepsilon^{\prime \prime}"]
& \text{}
 \\
 & F \arrow[r,equal] \arrow[d,dashed,"\varepsilon^{\prime}"] & F \arrow[d,dashed,"\mathbf{E}^{1}(F{,}f^{\prime})(\varepsilon^{\prime})"]
 \\ & \text{}  & \text{}
\end{tikzcd}
\end{center}
in $(\mathscr{C},\mathbf{E}^{1},\mathfrak{r})$ and an $\mathbf{E}^{1}$-extension $\varepsilon^{\prime \prime} \in \mathbf{E}^{1}(E,A)$ realised by the sequence $A \overset{h}{\longrightarrow} C \overset{h^{\prime}}{\longrightarrow} E$, satisfying the following compatibilities
\begin{enumerate}[(i)]
\item $\mathfrak{r}(\mathbf{E}^{1}(F,f^{\prime})(\varepsilon^{\prime}))  = [D \overset{d}{\longrightarrow} E \overset{e}{\longrightarrow} F].$

\item $\mathbf{E}^{1}(d,A)(\varepsilon^{\prime \prime}) = \varepsilon.$

\item $\mathbf{E}^{1}(E,f)(\varepsilon^{\prime \prime}) = \mathbf{E}^{1}(e,B)(\varepsilon^{\prime}).$
\end{enumerate}
This shows that the triple $(\mathscr{C},\mathbf{E}^{1},\mathfrak{r})$ satisfies (ET4). The proof showing that $(\mathscr{C},\mathbf{E}^{1},\mathfrak{r})$ satisfies (ET4)$^{\text{op}}$ is dual. 
\end{proof}
\end{proposition}

\begin{lemma}\label{NPTriang}\cite[Proposition 3.22(b)]{NakaokaPalu}
Let $\mathscr{A}$ be an additive category with an auto-equivalence $[1]$, and set $\mathbb{F}(-,-) = \mathscr{A}(-,-[1])$. If we are given an $\mathbb{F}$-triangulation $\mathfrak{t}$ of $\mathscr{A}$. Define that
\begin{center}
\begin{tikzcd}
A \arrow[r,"x"]
& B \arrow[r,"y"]
& C \arrow[r,"\d"]
& A[1]
\end{tikzcd}
\end{center}
is a distinguished triangle if and only if
$\mathfrak{t}(\d) = [A \overset{x}{\longrightarrow} B \overset{y}{\longrightarrow} C].$ Denote this class of distinguished triangles by $\Delta$. Then $(\mathscr{C},[1],\Delta)$ is a triangulated category. 
\end{lemma}

\textbf{Proof of the backward direction of Theorem \ref{mainTheorem}}
\begin{proof} Let $(\mathscr{C},\mathbb{E},\mathfrak{r})$ be an extriangulated category where for every object $X \in \mathscr{C}$ the morphism $X \rightarrow 0$ is an $\mathbb{E}$-inflation and the morphism $0 \rightarrow X$ is an $\mathbb{E}$-deflation. By Proposition \ref{SigmaAutoEquiv}, there is an auto-equivalence $\Sigma \colon \mathscr{C} \rightarrow \mathscr{C}$ and by Lemma \ref{ET1prf}, $\mathbf{E}^{1}(-,-) := \mathscr{C}(-,\Sigma -)$ is a biadditive functor. By Proposition \ref{ET2prf}, there is a correspondence $\mathfrak{r}$ which associates an equivalence class $\mathfrak{r}(\varepsilon)$ to any extension $\varepsilon \in \mathbf{E}^{1}(C,A)$ for any objects $A,C \in \mathscr{C}$. Moreover, the correspondence $\mathfrak{r}$ is an additive realisation. By Proposition \ref{ET3prf} the triple $(\mathscr{C},\mathbf{E}^{1},\mathfrak{r})$ satisfies (ET3) and (ET3)$^{\text{op}}$. By Proposition \ref{ET4prf} the triple $(\mathscr{C},\mathbf{E}^{1},\mathfrak{r})$ satisfies (ET4) and (ET4)$^{\text{op}}$. To summarise, we have an $\mathbf{E}^{1}$-triangulation of $\mathfrak{r}$ of $\mathscr{C}$. So by Lemma \ref{NPTriang}, $\mathscr{C}$ has a structure of a triangulated category $(\mathscr{C},\Sigma, \Delta)$, where $\Delta$ is the set of distinguished triangles as defined in Lemma \ref{NPTriang}. By Lemma \ref{NPTriang}, $\Delta$ is $\mathbf{E}^{1}$-compatible, by Lemma \ref{lemmaE1}, 
$\mathbf{E}^{1}(X,Y) \cong \mathbb{E}(X,Y)$, so $\Delta$ is $\mathbb{E}$-compatible. 

This completes the proof of Theorem \ref{mainTheorem}.
\end{proof}

\bibliographystyle{plain}
	\bibliography{IdempotentCompletion}
\end{document}